\documentclass[]{article}

\title{Convergence and generalization of a recursion equation for primes}
\author{James Haley}
\usepackage{amsmath}
\usepackage{amsfonts}
\usepackage{amssymb}
\usepackage{graphicx}

\begin{document}

\newtheorem{thm}{Theorem}

\maketitle

\section{Introduction}

There is no known ``nice'' function which produces the prime numbers. However, there are formulas from which, given the first few primes, one can obtain the remaining primes. In a paper by Keller [1], the following recursion equation is given:
\begin{equation*}
 p_{n+1}= \lim\limits_{s \to +\infty} \left[ \sum\limits_{j=1}^{2p_{n}-1} j^{-s} - \prod\limits_{k=1}^{n}(1-p_{k}^{-s})^{-1}\right] ^{-1/s}
\end{equation*}
where $p_{1}, p_{2}, \ldots , p_{n}$ are the first $n$ primes. This formula is proven using the Riemann zeta function. Our primary result in this paper is the demonstration that a similar formula is true for Dirichlet $L$-functions. To state our result precisely, we need to define a Dirichlet character. 

A Dirichlet character $\chi :\mathbb{Z}\rightarrow \mathbb{C}$ is a multiplicative function that has an associated modulus, say $k$. The function is periodic with $\chi(n) =\chi(n+mk)$ for any integer $m$, and $\chi (n) = 0$ if and only if $gcd(n,k)>1$. An example of a Dirichlet character is given in Fig. 1. 

{\scriptsize \textbf{Fig. 1:} The Dirichlet characters modulo 5. }

\begin{center}
\begin{tabular}{|c|c|c|c|c|c|}
\hline
\textbf{$\chi (n)$} & 0 & 1 & 2 & 3 & 4 \\ \hline
$\chi _{1}(n)$ & 0 & 1 & 1 & 1 & 1 \\ \hline
$\chi _{2}(n)$ & 0 & 1 & $i$ & $-i$ & -1 \\ \hline
$\chi _{3}(n)$ & 0 & 1 & -1 & -1 & 1 \\ \hline
$\chi _{4}(n)$ & 0 & 1 & $-i$ & $i$ & -1 \\ 
\hline
\end{tabular}
\end{center}

\vspace{3 pt}

\begin{thm}
Let $p_{1}, p_{2}, \ldots , p_{n}$ be the first $n$ prime numbers and let $\chi$ be a Dirichlet character. Then
\begin{equation}
p_{n+1} = \lim\limits_{s\rightarrow \infty} \left| \sum\limits_{j=1}^{2p_{n}-1} \dfrac{\chi (j)}{j^{s}} - \prod\limits_{k=1}^{n}(1-\frac{\chi (p_{k})}{p_{k}^{s}})^{-1}\right| ^{-1/s}.
\end{equation}
\end{thm}
Note that Keller's result is the case when $\chi$ is the trivial character with $\chi (n)=1$ for all $n$. This suggests a comparison of the convergence of various $\chi$ to see which converge to $p_{n+1}$ faster. We do this in Section 4 by examining the error function:
\begin{multline*} \begin{split} D_{n}(s,\chi ) :=& \left| p_{n+1} - \left[ \sum\limits_{j=1}^{2p_{n}-1} j^{-s} - \prod\limits_{k=1}^{n}(1-p_{k}^{-s})^{-1}\right] ^{-1/s}\right| \\ &- \left| p_{n+1} - \left| \sum\limits_{j=1}^{2p_{n}-1}\dfrac{\chi (j)}{j^{s}} - \prod\limits_{k=1}^{n}(1-\frac{\chi (p_{k})}{p_{k}^{s}})^{-1}\right| ^{-1/s} \right| .\end{split} \end{multline*}
Our numerical experiments indicate that certain $\chi$ cause the formula to converge faster, while the trivial character converges faster than other $\chi$.

In Section 3 we observe the convergence of Keller's formula by studying the error function
\begin{equation*} E_{n}(s):= \left| p_{n+1} - \left[ \sum\limits_{j=1}^{2p_{n}-1} j^{-s} - \prod\limits_{k=1}^{n}(1-p_{k}^{-s})^{-1}\right] ^{-1/s}\right| \end{equation*}
and find that 
\begin{equation*}
-\log E_{n}(s) \approx as+b
\end{equation*}
for some $a$ and $b$. A similar observation is made when we use a nontrivial $\chi$ and examine $E_{n}(s,\chi)$.

\section{Generalization of Keller's equation}
In this section we will prove Theorem 1. First note that 
\begin{multline} \prod\limits_{k=1}^{n}(1-\frac{\chi (p_{k})}{p_{k}^{s}})\cdot \sum\limits_{j=1}^{\infty}\dfrac{\chi (j)}{j^{s}} = \prod\limits_{k=1}^{n}(1-\frac{\chi (p_{k})}{p_{k}^{s}})\cdot \prod\limits_{j=1}^{\infty}(1-\frac{\chi (p_{j})}{p_{j}^{s}})^{-1}\\ = \prod\limits_{j=n+1}^{\infty}(1-\frac{\chi (p_{j})}{p_{j}^{s}})^{-1} = \ldots = \prod\limits_{j=n+1}^{\infty}(1+\frac{\chi (p_{j})}{p_{j}^{s}}+\frac{\chi (p_{j})^{2}}{p_{j}^{2s}}+\ldots )\end{multline}
 The final product in (2) is asymptotic to $1+\frac{\chi (p_{n+1})}{(p_{n+1})^{s}}$ as $s\rightarrow \infty$, since its expansion is the sum of $1+\frac{\chi (p_{n+1})}{(p_{n+1})^{s}}$ and a series of terms whose denominators are primes raised to powers that are increasing multiples of $s$. So we have\begin{equation*}
  \lim\limits_{s\rightarrow \infty} \prod\limits_{k=1}^{n}(1-\frac{\chi (p_{k})}{p_{k}^{s}})\cdot \sum\limits_{j=1}^{\infty}\dfrac{\chi (j)}{j^{s}} = \lim\limits_{s\rightarrow \infty} 1+\frac{\chi (p_{n+1})}{(p_{n+1})^{s}} \end{equation*}
  which implies that   
  \begin{equation} p_{n+1}^{-s} \sim \left[ \prod\limits_{k=1}^{n}(1-\frac{\chi (p_{k})}{p_{k}^{s}})\sum\limits_{j=1}^{\infty}\dfrac{\chi (j)}{j^{s}} -1\right] \left[ \chi (p_{n+1})\right] ^{-1}\end{equation}
  as $s\rightarrow \infty$. (Note that $\chi (p_{n+1})=0$ if and only if $p_{n+1}$ divides the modulus of $\chi$, so one can choose a modulus such that there is no division by zero.) Then, taking the $-1/s$ power of (3), we obtain
  \begin{multline} p_{n+1} = \lim\limits_{s\rightarrow \infty} p_{n+1} = \lim\limits_{s\rightarrow \infty} \left| \prod\limits_{k=1}^{n}(1-\frac{\chi (p_{k})}{p_{k}^{s}})\sum\limits_{j=1}^{\infty}\dfrac{\chi (j)}{j^{s}} -1\right| ^{-1/s} \left| \chi (p_{n+1})\right| ^{1/s}.\end{multline}
   But $|\chi (p_{n+1})|^{1/s} = (1)^{1/s} = 1$, so (4) yields
   \begin{equation}
   p_{n+1} = \lim\limits_{s\rightarrow \infty}\left| \prod\limits_{k=1}^{n}(1-\frac{\chi (p_{k})}{p_{k}^{s}})\sum\limits_{j=1}^{\infty}\dfrac{\chi (j)}{j^{s}} -1\right| ^{-1/s}.
   \end{equation}

\medskip
 
We can factor (5) to obtain
 \begin{multline}\left| \prod\limits_{k=1}^{n}(1-\frac{\chi (p_{k})}{p_{k}^{s}})\sum\limits_{j=1}^{\infty}\dfrac{\chi (j)}{j^{s}} -1\right| ^{-1/s}\\ = \left| \sum\limits_{j=1}^{\infty}\dfrac{\chi (j)}{j^{s}} - \prod\limits_{k=1}^{n}(1-\frac{\chi (p_{k})}{p_{k}^{s}})^{-1}\right| ^{-1/s} \left| \prod\limits_{k=1}^{n}(1-\frac{\chi (p_{k})}{p_{k}^{s}})\right| ^{-1/s}\end{multline}
 and $|\prod\limits_{k=1}^{n}(1-\frac{\chi (p_{j})}{p_{j}^{s}})|^{-1/s}$ is asymptotic to 1 as $s\rightarrow \infty$. Then (6) implies
 \begin{equation}
 p_{n+1} =  \lim\limits_{s\rightarrow \infty} \left| \sum\limits_{j=1}^{\infty}\dfrac{\chi (j)}{j^{s}} - \prod\limits_{k=1}^{n}(1-\frac{\chi (p_{k})}{p_{k}^{s}})^{-1}\right| ^{-1/s}.
 \end{equation}
 
\medskip
 
The sum in (7) can be made finite to obtain Theorem 1. Consider
\begin{align} & \sum\limits_{j=1}^{2p_{n}-1} \dfrac{\chi (j)}{j^{s}} - \prod\limits_{k=1}^{n}(1-\frac{\chi (p_{k})}{p_{k}^{s}})^{-1}=\cdots \\
& = \sum\limits_{j=1}^{2p_{n}-1} \dfrac{\chi (j)}{j^{s}} - \prod\limits_{k=1}^{n}(1+ \frac{\chi (p_{k})}{p_{k}^{s}}+ \frac{\chi (p_{k}^{2})}{p_{k}^{2s}}+ \frac{\chi (p_{k}^{3})}{p_{k}^{3s}}+ \cdots )\\
& = \sum\limits_{\substack{j>1\\ p_{i}\nmid j}}^{2p_{n}-1}\dfrac{\chi (j)}{j^{s}} - \sum\limits_{M}^{}\dfrac{c_{M}}{M^{s}}\end{align}
 for $i=1,2,...,n$, where each $M$ is a multiple of some $p_{i}$ with $M\geq 2p_{n}$ and $c_{M} = \prod\limits_{a,i\ s.t.\ p_{i}^{a}\parallel M}^{}\chi (p_{i}^{a})$. In the final difference (10), the smallest term of the first sum is $p_{n+1}^{-s}\chi (p_{n+1})$, and the smallest term of the second sum is $(2p_{n})^{-s}\chi (2)\chi (p_{n})$. Since $p_{n+1} < 2p_{n}$, the difference will converge to $p_{n+1}^{-s}$ as $s\rightarrow \infty$, and Theorem 1 follows.
 
\section{Convergence of Keller's equation}

From Keller [1], we have the recursion equation for primes
\begin{equation*} p_{n+1}= \lim\limits_{s \to +\infty} \left[ \sum\limits_{j=1}^{2p_{n}-1} j^{-s} - \prod\limits_{k=1}^{n}(1-p_{k}^{-s})^{-1}\right] ^{-1/s}\end{equation*}
where $p_{1}, p_{2},..., p_{n}$ are the first $n$ prime numbers. A natural question that arises from this formula is how quickly does it converge? For a given $s$ and a given $n$, consider the error value
\begin{equation*} E_{n}(s):= \left| p_{n+1} - \left[ \sum\limits_{j=1}^{2p_{n}-1} j^{-s} - \prod\limits_{k=1}^{n}(1-p_{k}^{-s})^{-1}\right] ^{-1/s}\right| \end{equation*}
Keller's result says that $E_{n}(s)\rightarrow 0$ as $s\rightarrow \infty$. The error values decrease at a highly regular rate, exhibiting a negative exponential pattern. Plotting $-\log (E_{n=2}(s))$, with $\log$ denoting the natural logarithm, for $20\leq s\leq 500$, we obtain a series of points that appear to be	 arranged linearly (see Fig. 2). The slope of the best-fit line of these points was approximately .2, indicating that each successive $s$ yields an estimate roughly 20\% closer to the actual prime than the previous $s$. Thus for each fixed $n$, it seems that, $\forall s\gg 0$, 
\begin{equation*} -\log E_{n}(s) \approx as+b \end{equation*} 
for some constants $a$ and $b$ which depend on $n$. In later work, we expect to determine the constants $a$ and $b$.

{\scriptsize \textbf{Fig. 2:} Plotting $-log(E_{2}(s))$ with $s$ on the horizontal axis.}

\begin{center}
\includegraphics[height=96 pt]{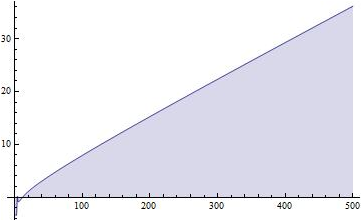}
\end{center}

The linear model was not perfect. There were some anomalies for $s<20$ that seem to be the result of taking the negative logarithm of a large error value. Additionally, data was only collected for $1\leq s\leq 500$, as larger $s$ took too long to compute. There was no indication that $s>500$ would produce different results. 

The pattern appeared for $n=2,3,...,20$. For each such $n$ and $s\in \mathbb{Z}$,\\ $1\leq s\leq 150$, $-\log E_{n}(s)$ had a correlation coefficent greater than .99. When the slopes of the best fit lines were plotted, they formed a vaguely negative exponential band (see Fig. 3). This indicates that, for larger primes, each successive $s$-value improves the accuracy of the Keller limit less than each successive\\ $s$-value does for smaller primes. This is somewhat expected. However, there is some oscillation within the band, so this does not hold absolutely. 

{\scriptsize \textbf{Fig. 3:} Plotting the slopes of $–log(E_n(20))$ with n along the horizontal axis.}

\begin{center}
\includegraphics[height=96 pt]{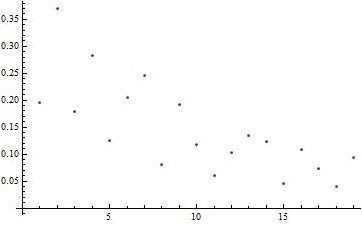}
\end{center}

The same linearity was found when we took the negative logarithm of
\begin{equation*}
E_{n}(s;\chi) := \left| p_{n+1} - \left| \sum\limits_{j=1}^{2p_{n}-1}\dfrac{\chi (j)}{j^{s}} - \prod\limits_{k=1}^{n}(1-\frac{\chi (p_{k})}{p_{k}^{s}})^{-1}\right| ^{-1/s} \right|,
\end{equation*}
the error values for the new version of the equation. It should be noted that there were some disturbances in the linearity for $s>300$ for certain $\chi$, probably arising from the inability of the program to efficiently compute complex absolute values and $s^{th}$-roots for large $s$. 

\section{Comparing convergence of the two equations}

It is natural to also study the convergence of the different forms of the equations. We examined the error differences, where
\begin{multline*} \begin{split} D_{n}(s,\chi ) :=& \left| p_{n+1} - \left[ \sum\limits_{j=1}^{2p_{n}-1} j^{-s} - \prod\limits_{k=1}^{n}(1-p_{k}^{-s})^{-1}\right] ^{-1/s}\right| \\ &- \left| p_{n+1} - \left| \sum\limits_{j=1}^{2p_{n}-1}\dfrac{\chi (j)}{j^{s}} - \prod\limits_{k=1}^{n}(1-\frac{\chi (p_{k})}{p_{k}^{s}})^{-1}\right| ^{-1/s} \right|. \end{split} \end{multline*}
 Fixing $s=50$ (somewhat arbitrarily - large enough to produce accurate results but small enough for fast computations), we evaluated $D_{n}(s,\chi )$ for $n = 3,4,5,6,7,8$ and for all $\chi$ modulo 4, 5, or 8 and those $\chi$ modulo 9 with complex values (see Fig. 4). All error differences were very small, indicating that both equations converge similarly. The error differences were both positive and negative, so neither equation consistently converges faster. Some $\chi$ yielded the same error differences as other $\chi$. In most cases, this was the result of each such $\chi$ being equivalent. There were some $\chi$, such as $\chi _{2}$ and $\chi _{4}$ modulo 5, that were different functions, but yielded the same error differences. However, their outputs had the same complex absolute value, which would explain why the difference errors were the same. With modulus 4,5, and 8, the $D_{3}(50,\chi )$ were equal for each $\chi$, and $D_{3}(50,\chi )$ were half that value for each $\chi$ tested with modulus 9. It is unclear why this occurred. 

{\scriptsize \textbf{Fig. 4:} $E_n(50, \chi )$ for $n = 3, 4, 5, 6, 7, 8$ and $\chi$ with modulo 4, 5, 8, and 9 }

\centerline{
\begin{tabular}{|c|c|c|c|c|c|c|}
\hline
\textbf{Prime:} & $3^{rd}$ & $4^{th}$ & $5^{th}$ & $6^{th}$ & $7^{th}$ & $8^{th}$ \\ \hline
$\chi _{1} \textit{mod 4}$ & $2.518\cdot 10^-9$ & $-0.000001277$ & $-9.921E-13$ & $-2.063E-10$ & $-9.287E-14$ & $-6.239E-12$ \\ \hline
$\chi _{2} \textit{mod 4}$ & $2.518\cdot 10^-9$ &$ -0.00000137$ & $2.988\cdot 10^-9 $&$ -0.000005125 $& $4.994\cdot 10^-10$ &$ 3.034\cdot 10^-7$ \\ \hline
$\chi _{1} \textit{mod 5}$ & $2.518\cdot 10^-9$ &$ -4.049E-08 $& $-1.641E-15 $&$ -1.572E-13 $&$ -2.063E-14 $& $-4.44E-13$ \\ \hline
$\chi _{2} \textit{mod 5}$ & $2.518\cdot 10^-9$ &$ 0.00004926 $&$ 1.494\cdot 10^-9 $&$ 0.001302 $&$ 0.00002698 $&$ 4.410\cdot 10^-6$ \\ \hline
$\chi _{3} \textit{mod 5}$ & $2.518\cdot 10^-9$ &$ -0.000002607$ &$ 2.989\cdot 10^-9$ &$ -0.000004939$ &$ -1.437E-09 $&$ -1.15E-11$ \\ \hline
$\chi _{4} \textit{mod 5}$ & $2.518\cdot 10^-9$ &$ 0.00004926 $& $1.494\cdot 10^-9 $&$ 0.001302 $& $0.00002698 $& $4.410\cdot 10^-6$ \\ \hline
$\chi _{1} \textit{mod 8}$ & $2.518\cdot 10^-9$ &$ -0.000001277 $& $-9.921E-13 $&$ -2.063E-10$ &$ -9.287E-14 $& $-6.239E-12$ \\ \hline
$\chi _{2} \textit{mod 8}$ & $2.518\cdot 10^-9$ &$ -0.000001289 $& $-1.59E-12 $& $1.847\cdot 10^-7 $&$ -1.954E-09 $& $-6.239E-12$ \\ \hline
$\chi _{3} \textit{mod 8}$ & $2.518\cdot 10^-9$ &$ -0.00000137 $& $2.988\cdot 10^-9 $&$ -0.000005125$ & $4.994\cdot 10^-10 $&$ 3.034\cdot 10^-7$ \\ \hline
$\chi _{4} \textit{mod 8}$ & $2.518\cdot 10^-9$ &$ -0.000001358$ &$ 2.987\cdot 10^-9 $&$ -0.000004939$ &$ -1.455E-09 $&$ 3.034\cdot 10^-7$ \\ \hline
$\chi _{2} \textit{mod 9}$ & $1.259\cdot 10^-9$ &$ 0.00002397 $& $1.966\cdot 10^-7 $& $0.\cdot 10^-5 $& $0.\cdot 10^-3 $& \\ \hline
$\chi _{3} \textit{mod 9}$ & $1.259\cdot 10^-9$ &$ 0.00002525 $&$ 1.951\cdot 10^-7 $& $0.\cdot 10^-5 $& $0.00001349 $& $2.433\cdot 10^-6$ \\ \hline
$\chi _{5} \textit{mod 9}$ & $1.259\cdot 10^-9$ &$ 0.00002525 $&$ 1.951\cdot 10^-7 $& $0.\cdot 10^-5 $& $0.00001349 $&$ 2.433\cdot 10^-6 $\\ \hline
$\chi _{6} \textit{mod 9}$ & $1.259\cdot 10^-9$ &$ 0.00002397 $&$ 1.966\cdot 10^-7 $& $0.\cdot 10^-5 $&$ 0.\cdot 10^-3 $& \\ 
\hline
\end{tabular}
}

\bigskip

Acknowledgments: I would like to thank Prof. C. Douglas Haessig for his guidance throughout this project. I would also like to thank Prof. Steven Gonek for his thoughts on the linearity of the error functions and the reviewers for their suggestions. 

\section*{References}
[1] Keller, Joseph B. ``A recursion equation for prime numbers.'' arXiv: 0711.3940v2. 5 October 2008.

\end{document}